\documentclass[english,12pt]{article}
\usepackage[totalwidth=430pt,totalheight=595pt]{geometry}
\usepackage{babel,amsmath,amssymb,amsbsy,amsfonts,latexsym,amsthm,epsf,array,colordvi,xy,fancyhdr,mathrsfs}
\input xy
\xyoption{all}

\def\Y#1{}

\newtheorem{teor}{Theorem}[section]
\newtheorem{lemma}[teor]{Lemma}
\newtheorem{prop}[teor]{Proposition}

\newcommand{\ac}{\`{a} }

\newcommand{\graor}{\widetilde{\mathbb{G}\mathrm{r}}_{3}}
\newcommand{\grasc}[1]{\mathbb{G}\mathrm{r}_{2}(\C^{\,#1})}
\newcommand{\grasr}[1]{\widetilde{\mathbb{G}\mathrm{r}}_{4}(\R^{\,#1})}
\newcommand{\la}{\mathfrak{g}}

\newcommand{\sud}{\mathfrak{su}(2)}
\newcommand{\sut}{\mathfrak{su}(3)}
\newcommand{\sot}{\mathfrak{so}(3)}

\newcommand{\ta}{\mathfrak{t}}
\newcommand{\C}{\mathbb{C}}
\newcommand{\Ha}{\mathbb{H}}
\newcommand{\R}{\mathbb{R}}
\newcommand{\Z}{\mathbb{Z}}
\newcommand{\N}{\mathbb{N}}
\newcommand{\Pro}{\mathbb{P}}

\newcommand{\si}{\Sigma}

\newcommand{\cons}{\mathbf{c}}
\newcommand{\princ}{\mathcal{P}}

\newcommand{\grad}{\mathrm{grad}}

\newcommand{\vspan}{\mathrm{span}}

\newcommand{\ds}{\displaystyle}
\newcommand{\ts}{\textstyle}
\newcommand{\double}{\mathscr{D}}
\newcommand{\sM}{\mathscr{M}}
\newcommand{\flag}{\mathbb{F}}
\newcommand{\ssu}{\mathbb{L}}
\newcommand{\aw}{\mathbb{A}}
\newcommand{\PP}{\mathbb{P}}
\renewcommand{\SS}{\mathbb{S}}

\newcommand{\Tr}{\mathrm{Tr}}

\newcommand{\uuu}{\mathfrak{u}}
\newcommand{\uo}[1]{\mathfrak{u}\lower1.5pt\hbox{$_{#1}$}}

\newcommand{\torus}{\mathrm{T}^{2}}
\newcommand{\weyl}{W}
\newcommand{\slice}{V}
\newcommand{\ugen}{\mathbf{u}}
\newcommand{\vgen}{\mathbf{v}}

\newcommand{\weight}{\mathbf{z}}

\newcommand{\xxbar}{\Gamma}
\newcommand{\hypers}{\mathcal{H}}

\begin{document} 
\title{\textsc{\bfseries Eight-dimensional SU(3)-manifolds\\ 
of cohomogeneity one }}
\author{\textsc{Andrea Gambioli}}
\date{}
\maketitle
\pagestyle{plain}

\abstract{\footnotesize In this paper, we classify $8$-dimensional
  manifolds $M$ admitting an $SU(3)$ action of cohomogeneity one such
  that (i) $M$ is simply connected and the orbit space $M/G$ is
  isomorphic to $[0,1]$, and (ii) $M/G\cong S^{1}$ and the principal
  orbits are simply connected. We discuss applications to the study of
  the group manifold $SU(3)$ and to $8$-dimensional quaternion-K\"ahler
  spaces.}\vspace{4mm}


\section{Introduction}\label{introduction}

Let $M$ be a differentiable manifold, and $G$ a compact semisimple
group acting smoothly on $M$. Then $M$ is said to be a
\emph{cohomogeneity-one} $G$-space if the principal orbits are
codimension-one submanifolds. A result due to Mostert \cite{mostert}
asserts that the quotient space $M/G$ is isomorphic to $[0,1]$ or to
$S^{1}$ if $M$ is compact, to $[0,1)$ or $\R$ if $M$ is non-compact.
In the case of the interval $[0,1]$, there are precisely two singular
orbits corresponding to the endpoints.

Manifolds with a cohomogeneity-one group action have been increasingly
studied in recent years. This is mainly due to the fact that many
problems concerning the existence of $G$-invariant structures on them
can be reduced to ODE's, which are sometimes straightforward to
handle. As typical examples, we cite \cite{beber}, \cite{brysal},
\cite{dancwang1}, in which such techniques were used to construct
Einstein metrics and examples of metrics with exceptional holonomy.

More recently, cohomogeneity-one qua\-ter\-nion-K\"{a}hler and
hyperk\"{a}hler manifolds were classified in \cite{dancswan1},
\cite{dancswan2}. General criteria for the classification of
cohomogeneity-one manifolds were also developed in \cite{alek-alek},
\cite{alek-alek2}, and used to partially classify such manifolds with $\chi(M)>0$
(and a corresponding family of quaternion-K\"{a}hler manifolds) in \cite{alek-pod}.
Cohomogeity-one $SU(3)$ manifolds of dimension 7 are the subject of
\cite{podestaverdiani}.

In this paper, we shall focus on $8$-dimensional simply-connected
smooth manifolds admitting an action of $SU(3)$ of cohomogeneity one.
The interest in these manifolds arises from the following
considerations. Firstly, the 8-dimensional quaternion-K\"{a}hler (QK) spaces
\begin{equation}\label{wolfspaces}
\Ha\Pro^{2},\qquad \grasc{4}\cong\grasr{6},\qquad G_{2}/SO(4),
\end{equation} remarkably all admit an $SU(3)$-action of cohomogeneity
one. 
(See \cite{wolf65}, \cite{alek3}, \cite{poonsal1} for the theory of
such Wolf spaces.) 
In \cite{gambioli1}, the author studied the moment mapping $\mu$ of a QK space
into the Grassmannian $\graor(\la)$ of oriented 3-planes in the Lie
algebra $\la$. Whilst $\mu$ is a branched covering of $G_2/SO(4)$ onto
its image in $\graor(\sut)$ \cite{kobswn1}, we shall point out  
that the first two spaces in (\ref{wolfspaces}) give rise to
7-dimensional images.

Another observation is that $8$ is precisely the dimension of the Lie
group $SU(3)$ itself, and it is natural to ask whether there is a
cohomogeneity-one action of $SU(3)$ on itself. Whilst the Adjoint
action has cohomogeneity-two, a positive answer to the question comes
from a modification called the \emph{$A$-twisted action} or
\emph{$\sigma$-action} (see \cite{conlon}, \cite{hptt1} and
\cite{kollross}). For the case of $SU(3)$, this coincides with the
more elementary \emph{consimilarity action}, studied independently in
the theory of matrices \cite{hornjohnson}. In any case,
the tangent space at a generic point of an $8$-dimensional Riemannian manifold with
an isometric $SU(3)$-action of cohomogeneity-one can be naturally
identified with the Lie algebra $\sut$.

 Such considerations suggest the importance of setting these four
examples in a wider context, in order to understand more deeply their
common structure.  Although $SU(3)$ does not admit a global QK
structure, we show that it has features in common with
(\ref{wolfspaces}) that allow it to be regarded as an ``honorary Wolf
space''. For example, $SU(3)$ minus a 5-sphere is
$SU(3)$-diffeomorphic to $G_2/SO(4)$ minus a complex projective plane
$\C\Pro^2$, and we explain that open dense sets of both $SU(3)$ and
$\Ha\Pro^{2}$ are the total spaces of $S^1$ bundles over the vector
bundle $\Lambda_-^2\C\Pro^2$. The manifold $SU(3)$ admits an invariant
hypercomplex structure \cite{joyc1}, and a $PSU(3)$ structure in the
sense of \cite{hitc1}. The theory also has links with $Spin(7)$
structures \cite{gukospartong}.

In the present article, we classify compact $8$-dimensional
differentiable manifolds $M$ admitting a cohomogeneity-one $SU(3)$
action such that the quotient space $M/SU(3)$ is $[0,1]$.  In this
case, the generic orbit has type $SU(3)/H$ where the connected
component at the identity is $S^{1}$, and there are precisely two
singular orbits $M_{1},\,M_{2}$ of type $SU(3)/K_{i}$, $i=1,2$,
satisfying the relations $SU(3)\supset K_{i}\supset H$ (we refer the
reader to \cite{bredon} for this basic theory). We also give a
partial classification of the case $M/G\cong S^{1}$, where in almost
all cases, $M$ turns out to be a product of $S^1$ with an
Aloff-Wallach space, which is the principal orbit.

The paper is organized as follows. In Section
\ref{preliminaryresults}, we describe our approach to the
classification, along with some results concerning connected principal stabilizers 
and the sphere-transitive representations of $U(2)$ and $\torus$; in
the latter case are also discussed non-connected pricipal stabilizers, which can 
appear only in presence of this type of singular stabilizer.
In Section \ref{theclassification}, we carry out the classification
distinguishing two possible situations: the case in which both
singular stabilizers are connected (Theorem \ref{classu3}), and that
in which at least one is not connected (Proposition \ref{propstab2}).
Moreover, we discuss the case in which $M/G\cong S^{1}$ and the
principal orbits are simply connected (Theorem \ref{mgs1}).

In Section \ref{examplesandapplications}, we shall identify some of
the manifolds obtained during the classification, and discuss more
extensively the consimilarity action of $SU(3)$ on itself.
Afterwards, in Section~\ref{quotientsbycirclesubgroups}, we apply
ideas behind the classification results to discuss the QK moment
mappings induced on $\Ha\Pro^{2}$ and $\grasc{4}$ under the action of
$SU(3)$, and relate these $8$-dimensional manifolds with examples 
of $7$-dimensional $SU(3)$-manifolds via circle actions.


\section{Preliminary results}\label{preliminaryresults}

In general, for arbitrary $G$-manifolds $M$ with orbit space isomorphic to $[0,1]$ there are two singular orbits $M_{1}$, $M_{2}$ 
and a normal (or slice) representation for each of them; let $\slice$ denote such representation  
at a point $x$ of a singular orbit $M_{i}$; then the bundle obtained as the twisted product
\begin{equation*}
G\times_{\scriptscriptstyle K_{i}} \slice
\end{equation*}
is $G$-equivariantly isomorphic to a tube around $M_{i}$. If we
consider the corresponding disk bundle $D_{i}$, we can describe $M$ as
\begin{equation*}
M=M_{\phi}=D_{1}\cup_{\phi} D_{2},
\end{equation*}
where 
\begin{equation}\label{gluing}
\phi:\partial D_{1}\xymatrix{\ar[r]&}\partial D_{2}
\end{equation}
is a $G$-equivariant diffeomorphism identifying the points of the two
boundaries. The latter are precisely the principal orbits:
$\partial D_{i}\cong G/H$, where $H$ is the principal stabilizer.

In \cite{uchida}, Uchida used this approach in order to classify
cohomology complex projective spaces with a cohomogeneity-one action.
We cite his useful sufficient conditions to decide if the manifolds
obtained using different maps $\phi$ are isomorphic as $G$-spaces (see
\cite[Lemma 5.3.1]{uchida}): let $\phi,\psi:\partial D_{1}\to
\partial D_{2}$ be $G$-equivariant maps as in (\ref{gluing}); then
$M_{\phi}$ and $M_{\psi}$ are $G$-equivariantly diffeomorphic if one
of the following conditions are satisfied:

\begin{enumerate}
\item
$\phi$ and $\psi$ are $G$-diffeotopic, or
\item
$\psi\circ \phi^{-1}$ can be extended to a $G$-equivariant
diffeomorphism of $D_{1}$ on itself, or
\item
$\phi\circ \psi^{-1}$ can be extended to a $G$-equivariant diffeomorphism of $D_{2}$ on itself.
\end{enumerate}
Our problem can therefore be reduced to classifying automorphisms of the generic orbit $SU(3)/U(1)$ up to these conditions.
 
One can obtain $G$-equivariant automorphisms of homogeneous spaces
$G/H$ as follows: let $a\in N(H)$; then the map $\phi^{a}$ given by
\begin{equation}\label{equivhom}
\phi^{a}(gH)=ga^{-1}H
\end{equation}
is well defined and commutes with the left multiplication for elements $g\in G$. It can be shown that \emph{all} 
$G$-equivariant automorphisms of $G/H$ have this form (see \cite[Chap
I, Th. 4.2]{bredon}); we have therefore the identification
\begin{equation*}
\mathrm{Aut}_G(G/H)\cong \frac{N(H)}{H}\,.
\end{equation*}\smallbreak

Let us discuss in some detail the case that two $SU(3)$-spaces
obtained by distinct gluing maps are isomorphic (as $SU(3)$-spaces).
In general if $M_{\phi}$ and $M_{\psi}$ are two such manifolds, then
an equivariant morphism $\Phi:M_{\phi}\to M_{\psi}$ would restrict on
the two tubular neighborhoods to a couple of equivariant morphisms
$\phi^{a}$ and $\phi^{b}$, as described in (\ref{equivhom}), which
make the following diagram commutative:
\begin{equation}\label{gluingdiag}
\xymatrix{G/K_{1} & G/H\ar[l]_{p_{1}} & G/H\ar[r]^{p_{2}}\ar[l]_{\psi}&G/K_{2}\\
G/K_{1}\ar[u]^{\phi^{a}}&G/H\ar[u]^{\phi^{a}}\ar[l]^{p_{1}}& G/H\ar[u]_{\phi^{b}}\ar[r]_{p_{2}}\ar[l]^{\phi}& G/K_{2}\ar[u]_{\phi^{b}}}
\end{equation}
In this case, $a\in N(H)\cap N(K_{1})$ and $b\in N(H)\cap N(K_{2})$.
\smallbreak

In general we cannot expect to have the same map $\phi^{a}$ in the
first two columns of the diagram (similarly for $\phi^{b}$);
instead, for instance, we will have $\phi^{a}$ and $\phi^{a'}$ repectively at
$G/K_{1}$ and at $G/H$.  Nevertheless, the map $\Phi$ is always
diffeotopic (through $SU(3)$-invariant maps) to a map $\Phi'$ for wihich
$a$ and $b$ are constant in the respective tubular neighborhoods.  The
homotopy between $\Phi$ and $\Phi'$ can be described as follows: the
map $\Phi$ is identified on each tubular neighborhood by a continuous
function $\epsilon:[0,\frac12]\to N(H)$, so that $\Phi=
\phi^{\epsilon(t)}$; we can define $$\eta(t,s):=\epsilon((1-s)t)\,,$$
and $\phi^{\eta(t,s)}$ is the required homotopy. We also observe that, for instance, 
$a\in N(H)\cap N(K_{1})$ in general, because the map $\epsilon$ is
continuous and $N(H)$ is a closed subgroup.  \vspace{2mm}

In the sequel, we shall use the following notation to identify the
most commonly used homogeneous spaces:

\begin{equation*}
\hspace{30pt}\begin{array}{ll}
\ds \SS:=\frac{SU(3)}{SU(2)}, &\hbox{the 5-sphere},\\[10pt]
\ds \PP:=\frac{SU(3)}{S(U(2)\times U(1))}, &\hbox{the complex projective
  plane }\C\Pro^2,\\[10pt]
\ds \ssu:=\frac{SU(3)}{SO(3)}, &\hbox{the set of special Lagrangian
  subspaces in $\C^{3}$},\\[10pt]
\ds \aw:= \frac{SU(3)}{U(1)}, &\hbox{any Aloff-Wallach type space},\\[10pt]
\ds \flag:=\frac{SU(3)}{T^2}, &\hbox{the complex flag manifold}. 
\end{array}
\end{equation*}

We shall actually use $\aw$ to stand for any homogeneous space of the
form $SU(3)/U(1)$, even though the terminology ``Aloff-Wallach'' usually
excludes one case (we shall be more precise in the next section).  The
Lagrangian interpretation of $\ssu$ can be found in
\cite{harveylawson}, and is important for making more explicit some of
our constructions, such as finding geodesics from one singular orbit
to another.


\subsection{Connected principal stabilizers}

Principal stabilizers $H$ in our case are $1$-dimensional subgroups of $SU(3)$,
such that $H^{0}=U(1)$. The case in which $H^{0}=H$ will be particularly significant, so 
we will dedicate the first part of this section to it.
We begin by defining circle subgroups of $SU(3)$. Let $k,l$ be
integers, and let $U_{k,l}$ denote the subgroup (isomorphic to $U(1)$)
of $SU(3)$ consisting of matrices
\[\begin{pmatrix}
e^{k\imath t} & 0 & 0\\
0 & e^{l \imath t} & 0 \\
0 & 0 & e^{-(k+l)\imath t}\\
\end{pmatrix}\,.
\]
We shall denote the coset space $SU(3)/U_{k,l}$ by $\aw_{k,l}$.  Since
$U_{k,l}$ is unchanged when any common factor of $k,l$ is removed, we
may assume that they are coprime. The space $\aw_{k,l}$ is called an
\emph{Aloff-Wallach space} provided $kl(k+l)\ne0$, since the pairs
equivalent to $(1,-1)$ are excluded for geometrical reasons (they do
not satisfy the conditions that guarantee the existence of homogeneous
positively-curved metrics, see \cite{aw}). In our analysis, the
subgroups $U_{1,-1}$ will however play important roles.

Denote the 1-dimensional subalgebra of $\sut$ corresponding to
$U_{k,l}$ by $\uo{k,l}$. We consider the pair of orthogonal
subalgebras $\uo{1,-1},\uo{1,1}$ generated by the respective
elements
\begin{equation}\label{princu1}
\ugen=
\begin{pmatrix}
\imath & 0 & 0\\
0 & -\imath  & 0 \\
0 & 0 & 0\\
\end{pmatrix}\quad\text{and}\quad
\vgen= 
\begin{pmatrix}
\imath & 0 & 0\\
0 & \imath  & 0 \\
0 & 0 & -2\imath\\
\end{pmatrix}
\end{equation}
that together span a Cartan subalgebra $\ta$.  It can be shown
that $\ugen$ is a regular element, so it belongs to a unique Cartan
subalgebra $\ta\subset \sut$, namely that consisting of diagonal
elements; if $\alpha,\,\beta,\alpha+\beta$ denote the roots in
$\ta_{\C}$, we have that $\ugen$ corresponds to $\alpha$, and we can 
identify $$\vspan\{\alpha\}=-\imath\, \uo{1,-1}\quad \vspan\{\beta\}= -\imath\, \uo{1,0}\quad \vspan\{\alpha+\beta\}= -\imath\,\uo{0,1}\,.$$ 
On the other hand, $\vgen$ is a singular element and is contained in three
independent Cartan subalgebras $\ta,\,\ta_{1},\,\ta_{2}$; the
1-dimensional orthogonal complements
$\vgen^{\perp},\,\vgen^{\perp}_{1},\,\vgen^{\perp}_{2}$ then span
the subalgebra $\sud$ corresponding to the root $\alpha$. 
Each root has an orthogonal singular hyperplane, which in our notation
correspond to $\uo{1,1}, \uo{-2,1}$ and $\uo{1,-2}$.\vspace{2mm}

The first step to obtain our classification is that of determining the possible gluing maps
between two principal orbits. In the case that the principal stabilizer is connected, this
corresponds to identify the group $N(U(1))$: this depends from the way $U(1)$ is 
immersed in $SU(3)$, up to conjugacy. The subgroups $U_{1,-1}$ and $U_{1,1}$ represent 
distinguished cases, in this sense. 

\begin{lemma}\label{normu1}
The normalizer of $U_{k,l}$ in $SU(3)$ is given by
\begin{equation*}
N(U_{k,l})=
\begin{cases}
\torus\cup \tau\torus & \text{if $(k,l)=(1,-1)$}\\
S(U(2)\times U(1)) & \text{if $(k,l)=(1,1)$}\\
\torus & \text{if $(k\pm l)\ne0$.}\\
\end{cases}
\end{equation*}
Here, $\tau$ denotes an element of $SU(3)$ such that $Ad_{\tau}$ is an element in the Weyl group $\weyl$.
\end{lemma}
\proof For the first case, let $g\in N(U_{1,-1})$; then we also have $g\in N(\torus)$, as otherwise 
\begin{equation*}
U_{1,-1}\subset g\torus g^{-1}\neq \torus 
\end{equation*}
which is impossible as $\ugen$ is a regular element. Hence
$N(U_{1,-1})\subset N(\torus)$. It is well known that
\begin{equation*}
\weyl:=\frac{N(\torus)}{\torus}\cong \mathfrak{S}_{3} 
\end{equation*}
is the group of permutations on $3$ elements; it acts on the Cartan
subalgebra $\ta$ by permuting the three roots $\alpha,\,\beta$ and
$\alpha+\beta$. The only elements fixing the subspace $t\cdot\alpha$
corresponding to $\ugen$ are reflections about the hyperplane
$\ugen^{\perp}$, sending $\ugen$ to $-\ugen$ and swapping $\beta $
and $\alpha+\beta$, which can be represented by the the action
$Ad_{\tau}$ with $\tau$ an appropriate element in $SU(3)$.

In the second case, an element $g\in N(U_{1,1})$ that preserves
$\uo{1,1}$ must also preserve the centralizer $C(U_{1,1})=S(U(2)\times U(1))\cong U(2)$, 
so that $N(U_{1,1})\subset N(U(2))=U(2)$; the reverse inclusion is obvious.

In the final case, we just have to note that roots and their
orthogonal complements are the only eigenspaces for the elements of
$\weyl$. Hence the other regular elements in $\ta$ are normalized only
by $\torus\cong N(\torus)^{0}$.\ $\blacksquare$ \vspace{3mm}

As a consequence, we obtain the required isomorphisms for the coset
spaces parametrizing $SU(3)$-equivariant automorphisms of principal
orbits. Firstly,
\begin{equation*}
\frac{N(U_{1,-1})}{U_{1,-1}}\cong U(1)\cup \tau U(1);
\end{equation*}
more explicitly, this group is generated by the matrices
\begin{equation}\label{normstabroot}
\begin{pmatrix}
e^{\imath t} & 0 & 0\\
0 & e^{\imath t} & 0 \\
0 & 0 & e^{-2 \imath t} \\
\end{pmatrix}
\quad\text{and}\quad
\begin{pmatrix}
0 & e^{\imath t} & 0\\
-e^{\imath t} & 0 & 0 \\
0 & 0 & e^{-2 \imath t} \\
\end{pmatrix}\,.
\end{equation}
For the second case,
\[\frac{N(U_{1,1})}{U_{1,1}}\cong SU(2),\]
and finally
\[\frac{N(U_{k,l})}{U_{k,l}}\cong U(1),\qquad (k\pm l)\ne0.\]

\noindent\emph{Remark.} We can already estimate the number of
$SU(3)$-equivariant diffeomorphism classes in some cases. For, if the
principal stabilizer is conjugate to $U_{1,-1}$ then, thanks to Lemma
\ref{normu1} and Uchida's condition 1, there are at most two such
classes (the number of connected components of
$N(U_{1,-1})/U_{1,-1}$). In the case of singular and all other regular
elements there is just one $SU(3)$-diffeomorphism
class.\vspace{1mm}

The next information that we need is knowledge of which tubular
neighborhoods can be built around a given singular orbit. To this aim,
we have to determine which representations of a singular stabilizer
are sphere-transitive, and associate to it the integers $k,l$
characterizing the corresponding principal stabilizer.


\subsection{$U(2)$ representations}

Let us concentrate now on the subgroup $S(U(2)\times U(1))\cong U(2)$
of $SU(3)$, classifying its sphere-transitive real
$4$-di\-men\-sio\-nal representations.

First we introduce some notation. Let $\si^n$ denote the complex
irreducible representation of $SU(2)$ on $\C^2$ of dimension $n+1$,
and $A^{m}$ the $U(1)$-re\-pre\-sen\-ta\-tion of weight $m$ with $m\in\Z$.

\begin{prop}\label{su2u1}
The real $4$-dimensional sphere-transitive representations $\slice$ of
$U(2)$ are given by
\begin{equation*}
\slice_{\C}\cong\C^{4}\cong\si^1\otimes (A^{m}\oplus A^{-m}),\quad m=2r+1,\,r\in\Z.
\end{equation*}
If $\{\ugen,\vgen\}$ is a basis for $\ta\subset \sud\oplus\uo{1,1}$
(see (\ref{princu1})), then the Lie algebra of the stabilizer of a
point $x\in S^{3}\subset \slice\cong\R^{4}$ has the form $(3\ugen+
m\vgen)^{\perp}$.
\end{prop}
\proof\, A consequence of the Peter-Weyl theorem is that a
 representation of $SU(2)\times U(1)$ necessarily has
the form
\begin{equation*}
\slice_{\C}\cong
\sum_{n,m}\Sigma^{n}\otimes A^{m}.
\end{equation*}
It is straightforward to see that the only possible case in which one
can obtain a sphere-transitive $4$-dimensional representation is given by
$\si^1\otimes (A^{m}+A^{-m})$. Now, the sums $1+m$ and $1-m$ must be
even in order to obtain an $S(U(2)\times U(1))\cong U(2)$
representation, as $SU(2)\times U(1)$ covers $U(2)$ in a two-to-one
manner. Hence $m$ must be odd.

Let us restrict the representation to the maximal torus $\torus$
contained in $U(2)$,
whose Lie algebra is $\ta=\vspan\{\ugen,\,\vgen\}$; then we obtain
\begin{equation*}
\slice_{\C}\cong (A^{1}+A^{-1})\otimes (A^{m}+A^{-m})\cong A^{m+1}+A^{-m-1}+A^{-m+1}+A^{m-1}.
\end{equation*}
The necessary real structure is effectively the tensor product $\jmath\otimes\jmath$
of the respective quaternionic structures on $\si^1$ and
$A^{m}+A^{-m}$. The latter act as the antilinear extensions of 
the maps $\jmath(x,y)=(y,-x)$ and $\jmath(e,f)=(f,-e)$ for $x,y$ a basis of $\si$ and $e,f$ a basis of $A^{m}+A^{-m}$;
the fixed point set is given by 
\begin{align*}
\slice=\>&\vspan\{x\otimes e+y\otimes f,\,x\otimes f-y\otimes e,\>\imath(x\otimes e-y\otimes f),\>\imath(x\otimes f+y\otimes e)\}\notag\\
=\>&\vspan\{w_{1},\,w_{2},\,w_{3},\,w_{4}\}.\notag
\end{align*}
Let us consider now the corresponding Lie algebra representation. We choose the point $w_{1}\in S^{3}$: 
then the Lie algebra $\sud\oplus \uo{1,1}$ acts on $w_{1}$ spanning the $3$-dimensional
tangent space of $S^{3}$. More explicitly, if
\begin{equation*}
\sud=\vspan\{v_{1},\,v_{2},\,v_{3}\}=\vspan\left\{\begin{pmatrix}\imath&0\\0&-\imath\end{pmatrix},\,\begin{pmatrix}0&\imath\\\imath&0\end{pmatrix}
,\,\begin{pmatrix}0&1\\-1&0\end{pmatrix}\right\}
\end{equation*}
then we obtain 
\begin{equation*}
v_{1}(w_{1})=w_{3},\quad v_{2}(w_{1})=w_{4},\quad v_{3}(w_{1})=w_{2}\,;
\end{equation*}
moreover the generator $\vgen\in\uo{1,1}$ acts by
\begin{equation*}
\vgen(w_{1})=m w_{3}\,.
\end{equation*} 
Returning to the inclusion $\sud\oplus\uo{1,1}\subset\sut$
and identifying $v_{1}=\ugen$, we can restrict to the Cartan 
subalgebra $\ta$; then the subspace spanned by 
$w_{1}, w_{3}$ in $\slice$ is an irreducible $\ta$-submodule, and the 
corresponding weight can be represented via the Killing metric by 
the vector
\begin{equation*}
\ts \mathbf{h}=\frac{1}{2}\ugen + \frac{m}{6}\vgen;
\end{equation*}
its kernel, which kills the vector $w_{1}$, is given by the
hyperplane $\mathbf{h}^{\perp}$: this is the Lie algebra of the 
stabilizer $U(1)$, hence the conclusion.\,$\blacksquare$
\vspace{2mm}

Let us analyse in more detail some examples for small $m$: for
$m=1$ we get the hyperplane $\uo{0,1}$ as the stabilizer's Lie
algebra; for $m=3$ we obtain a singular stabilizer $\uo{2,-1}$. 
For $m\geq5$ we get other generic
regular stabilizers, all belonging to the Weyl chambers 
delimited by $\uo{-2,1}$ and $\uo{1,-2}$; the limit stabilizing subalgebra 
for $m \to \infty$ corresponds to $\uo{1,-1}$. 
The same results are obtained for $m\leq0$, as the representation $A^{m}$
and $A^{-m}$ are isomorphic as real reprsentations.
\vspace{2mm}

In the sequel, let us use $\PP(m)$ to denote the bundle on
$\PP=\C\Pro^{2}$ obtained as the twisted product by the representation $\slice
=[\si^2\otimes (A^{m}+A^{-m})]$. Clearly $\PP(m)\cong \PP(-m)$, so we can restrict
to $m\in \N$.


\subsection{$\torus$ representations}

Let us discuss now the case of $\torus$ as a singular stabilizer:
we need to determine its sphere transitive $2$-dimensional representations in order
to classify the possible tubular neighborhoods around a singular orbit of type $\flag$. Let us choose for the standard Cartan subalgebra 
$\ta$ the basis formed by
\begin{equation}\label{basiscartan}
\ugen=
\begin{pmatrix}
\imath & 0 & 0\\
0 & -\imath  & 0 \\
0 & 0 & 0\\
\end{pmatrix}\quad\text{and}\quad
\ugen'= 
\begin{pmatrix}
0 & 0 & 0\\
0 & \imath  & 0 \\
0 & 0 & -\imath\\
\end{pmatrix}\,.
\end{equation}

Comparing this basis with that in (\ref{princu1}), we note that the
relation $\vgen=\ugen+2\ugen'$ holds, and that $\ugen,\ugen'$
correspond to the two roots $\alpha,\beta$; the parallelogram $P$
determined by $2\pi\ugen$ and $2\pi\ugen'$ is a fundamental domain for
the maximal torus $\torus$, which can therefore be described as
\begin{equation*}
\torus \cong \{\exp{s\ugen}\times\exp{t\ugen'}:s,t\in\R\}\,.
\end{equation*}

The $2$-dimensional spere-transitive real $\torus$-represenations $\slice$ are given by
\begin{equation*}
\slice\cong A^{p}\otimes A^{q} 
\end{equation*} 
for $p,q\in \Z$, with $(p,q)\neq(0,0)$ and $A^{p}\otimes A^{q}\cong A^{-p}\otimes A^{-q}$. Each of them is 
determined by a weight $\weight$ contained in $\ta$ such that 
\[
\langle\weight,\ugen\rangle=p,\qquad 
\langle\weight,\ugen'\rangle=q
\]
A basis for the integer lattice of such weights is given by 
\begin{equation*}\ts
\weight_{1}:=\frac{1}{3}\,\vgen\,,\quad 
\weight_{2}:=\frac{1}{3}\,(2\ugen+\ugen')\;\quad
\end{equation*} 
so that a generic weight has the form $\weight=p\,\weight_{1}+q\,\weight_{2}$ for $p,q\in\Z$, and the stabilizer
for the corresponding representation is given by $\weight^{\perp}$.
\vspace{3mm}

\noindent\emph{Observation.} The weights described in Proposition \ref{su2u1} are of this type: in fact
\begin{equation*}\ts
\frac{1}{2}\ugen + \frac{m}{6} \vgen=\frac{1}{2}\ugen + \frac{m}{6} (\ugen+2\ugen')=\frac{m-1}{2}\weight_{1}+\weight_{2}
\end{equation*}
for the choice $p=(m-1)/2$ and $q=1$ (recall that $l$ is odd): this is just the result of the reduction from $U(2)$
its maximal torus $\torus$. In fact the representation ring $R[U(2)]$
is isomorphic to the polynomial ring $\Z[\lambda_{1},\lambda_{2},\lambda_{2}^{-1}]$,
whereas $R[\torus]$ is isomorphic to
$\Z[\lambda_{1},\lambda_{1}^{-1},\lambda_{2},\lambda_{2}^{-1}]$; 
it is well known that the inclusion $\torus\subset U(2)$ induces
an injective map 
\begin{equation*}
\xymatrix{R[U(2)]\ar[r]&}R[\torus]
\end{equation*}
(see \cite{brocktd}). The image of this inclusion coincides with the
subring $$R[\torus]^{\weyl(U(2))}$$
of $\torus$ representations which
are invariant under the Weyl group $W(U(2))$; the latter is
isomorphic to $\Z_{2}\subset \weyl(SU(3))$, corresponding to a
reflection around one of the singular hyperplanes.  \vspace{3mm}

In this case, generic stabilizers might not be connected:

\begin{lemma}\label{lemmastabtorus}
For a $\torus$-representation of type $A^{p}\otimes A^{q}$ the generic stabilizer is
\begin{equation*}
U(1)\times \Z_{h}
\end{equation*} 
where $h=\gcd(p,q)$. 
\end{lemma}
\proof\, We can describe the representation by $(x,y)\mapsto
e^{2\pi\imath(px+qy)}$, where $(x,y)$ are coordinates with respect to (\ref{basiscartan}), after 
a suitable normalization, to be considered modulo $\Z^{2}$. The stabilizer is the solution of the
  equation
\begin{equation}\label{stabtorus}
px+qy=h,\qquad h\in \Z\,;
\end{equation}
so we have a $1$-dimensional solution for each $h$. On the other hand we can choose any $h_{1}$ and $h_{2}$ in $\Z$ 
so that $(x+h_{1},y+h_{2})$ is the same solution as $(x,y)$ on $\torus$, but for a different $h$. Therefore equation (\ref{stabtorus}) becomes
\begin{equation}\label{gcd}
px+qy=h-ph_{1}-qh_{2}\,.
\end{equation}
Let us suppose that $h=\gcd(p,q)>0$: then equation (\ref{gcd}) is equivalent to $px+qy=0$, as the $\gcd$ is precisely 
the smallest positive integer which can be obtained in the form $ph_{1}+qh_{2}$. Hence the solution $(x,y)$ for $h$ is also the 
solution for $0$; moreover this implies that if $0<h'<h$, then the solution $(x,y)$ for $h'$ is not a solution for $0$. This shows that the solutions are repeated 
modulo $h$, so that there are precisely $h$ distinct ones, each one isomorphic to a circle $U(1)$: altogether they form an
abelian subgroup, isomorphic to $U(1)\times \Z_{h}\subset \torus$. $\blacksquare$
\vspace{2mm}

We introduce some more notation at this point: we shall denote by
$\flag(p,q)$ a tubular neighborhood of a flag manifold obtained by a
slice representation $A^{p}\otimes A^{q}$ as explained above. We 
observe that $\flag(p,q)\cong \flag(-p,-q)$.

\def\Tone{
\begin{table}[h]
\begin{center}
\setlength{\extrarowheight}{1mm}
\setlength{\tabcolsep}{2mm}
\vspace{0.5cm}
\begin{tabular}{|c||c|c|c|c|c|}
\hline
   & $SU(2)$ & $U(2)$ &  $SO(3)$ & $\torus$\\[1pt]
\hline
\hline
  $\dim \slice$ & $3$ & $4$ & $3$ & $2$ \\[1pt]
\hline
  $\slice$ & $[\si^2]$ & $[\si\otimes(A^{l}\oplus A^{-l})]$ & $[\si^2]$ & $A^{p}\otimes A^{q}$ \\[1pt]
\hline
\end{tabular}
\parbox{320pt}{\caption{Connected singular stabilizers and
    corresponding slice representations}}
\end{center}
\label{slicerepr}
\end{table}
}

\def\Ttwo{
\begin{table}[h]
\begin{center}
\setlength{\extrarowheight}{1mm}
\setlength{\tabcolsep}{2mm}
\vspace{0.5cm}
\begin{tabular}{|c||c|c|c|c|}
\hline
  $_{M_{2} }\backslash^{M_{1}}$ & $\SS$ &  $\ssu$ & $\PP(l)$  \\[1pt]
\hline
\hline
  $\SS$ & $1$ & $1$  & $\delta_{l}^{1}$   \\[1pt]
\hline
 $\ssu$ & $1$ & $1$ & $ \delta_{l}^{1}$   \\[1pt]
\hline
$\PP(m)$ & $\delta_{m}^{1}$ & $ \delta_{m}^{1}$ & $\delta_{l}^{m}+\delta^{l}_{1}\delta_{m}^{1}$  \\[1pt]
\hline
\end{tabular}
\parbox{350pt}{\caption{Numbers of $SU(3)$-diffeomorphism classes of
    $8$-manifolds: 
singular stabilizers $SU(2)$, $U(2)$ and $SO(3)$ ($m,l$ odd)}}
\end{center}
\label{su3classes1}
\end{table}
}

\def\Tthree{
\begin{table}[h]
\begin{center}
\setlength{\extrarowheight}{1mm}
\setlength{\tabcolsep}{2mm}
\vspace{0.5cm}
\begin{tabular}{|c||c|c|c|c|}
\hline
  $_{M_{2} }\backslash^{M_{1}}$ & $\flag(l,m)$ & $\PP(l)$ & $\ssu$ & $\SS$  \\[1pt]
\hline
\hline
  $\flag(p,q)$ & $\delta_{l}^{p}\delta_{q}^{m}$ & $\delta^{p}_{(l-1)/2}\delta^{1}_{q}$ & 
$\delta^{p}_{0}\delta_{q}^{1}$ & $\delta^{p}_{0}\delta_{q}^{1}$\\[1pt]
\hline
\end{tabular}
\parbox{350pt}{\caption{Numbers of $SU(3)$-diffeomorphism
    classes of $8$-manifolds:
one singular stabilizer of type $\torus$ ($(p,q)\neq(0,0)$ and $\gcd(p,q)=1$)}}
\end{center}
\label{su3classes5}
\end{table}
}

\Tone

Regarding the singular stabilizers $SO(3)$ and $SU(2)$, we have a unique sphere-tran\-si\-tive
$3$-di\-men\-sio\-nal representation, namely the standard irreducible space $\R^3\cong
[\si^2]$. The complete list of possible slice representations for each connected
singular stabilizer is given in Table~1.


\section{The classification}\label{theclassification}

We are now in a position to present the main results of the paper,
classifying the possible ways of gluing together tubular
neighborhoods obtained from the singular orbits discussed in Section
\ref{introduction} and from the normal representations described in
Section \ref{preliminaryresults}.


\subsection{Connected singular stabilizers}

We focus first on the case that both the singular stabilizers
$K_1,K_2$ are connected. Connected subgroups of $SU(3)$ are in
one-to-one correspondence with Lie subalgebras of $\sut$. Note that
the two Lie subalgebras $\sot$ and $\sud\oplus\R$ are maximal
subalgebras. It is also well known that $\sot$ and $\sud$ are the
only $3$-dimensional subalgebras of $\sut$, up to conjugation.

Passing to subalgebras of $\sot $ and $\sud\oplus \R$, observe that
$\sud\cong\sot$ does not contain any subalgebra of dimension greater
than $1$; therefore we obtain only other two subalgebras, both
contained in $\sud\oplus\R$: namely $\sud$ and the Cartan subalgebra
$\ta$.

Let us also list here the normalizers of each corresponding connected
subgroup. Let $\Z_{3}$ denote the center of $SU(3)$, and (again)
$\weyl\cong\mathfrak{S}_3$ its Weyl group. Then
\[\begin{array}{c}
N(SU(2))=N(U(2))=U(2)\\[5pt]
N(SO(3))=SO(3)\times \Z_{3}\\[5pt] 
N(\torus)= \bigcup\limits_{\tau\in \weyl} \tau \torus
\end{array}\]\vspace{1mm}

\noindent\emph{Remark.} We shall not treat immediately the case of a
singular stabilizer $\torus$ with slice representation $A^{p}\otimes
A^{q}$ and $\gcd(p,q)\neq1$. In fact this can imply that the second
singular stabilizer is not connected even if $\torus$ is (because the
principal stabilizer turns out to be not connected, see Lemma
\ref{lemmastabtorus}), and this situation fits better in Subsection
\ref{nonconnectedstabilizers} (see Proposition \ref{propstab2}).
\vspace{3mm}

We can now state the main result of this section:

\begin{teor}\label{classu3}
  Tables~2 and 3 list respectively all the $SU(3)$-diffeomorphism
  classes of $8$-dimensional compact cohomogeneity-one
  $SU(3)$-manifolds
  with orbit space $[0,1]$ such that:\\[5pt]
  --- both stabilizers belong to the set $\{SU(2),\,U(2),\,SO(3)\}$,\\[5pt]
  --- one singular stabilizer is isomorphic to $\torus$ and the normal
  representation is $A^{p}\otimes A^{q}$ with $\gcd(p,q)=1$.
\end{teor}

\Ttwo

\Tthree
 
\vskip5pt

\proof Let us consider these connected singular stabilizers:
correspondingly we have a slice representation $\slice$ of dimension
$3,\,4,\,3,\,2$ (see Table~1); the representations involved must again
be of cohomogeneity one, or in other words the singular stabilizer
$K_{i}$ must act transitively on the unit sphere $S^{n-1}\subset
\slice$. Let us analyse the possibilities case by case.
\vspace{2mm}

 The cases of $SU(2)$ and $SO(3)$ are rather simple, as the only $3$-dimensional
representation of cohomogeneity one is the standard 3-dimensional
irreducible representation $\R^{3}\cong[\si^2]$, as already observed
at the end of Section \ref{preliminaryresults}; in this case the
principal stabilizer turns out to be one corresponding to $U_{1,-1}$. 
Therefore, thanks to Lemma \ref{normu1}, the normalizer is
$\torus \cup\tau \torus$ in both cases; on the other hand, it can be
shown that for both the singular orbits $\SS$ and $\ssu$, the component $\tau\torus$ of the normalizer $N(U(1))$ 
intesects $SU(2)\subset N(SU(2))$ and $SO(3)\subset N(SO(3))$ 
respectively (for instance in a point $x$ obtained by putting $t=\pi/2$
in an appropriate conjugate of the second element in
(\ref{normstabroot})). Hence any $SU(3)$ equivariant automorphism
of the principal orbit is diffeotopic to one which can be extended to an automorphism of the whole
tubular neighborhood (see Uchida's criteria in Section \ref{preliminaryresults}), so that
we have a unique $SU(3)$-equivariant diffeomorphism class of $M$
containing one of $\SS$ or $\ssu$ and another singular orbit $M_{2}$.\vspace{2mm}

 Let us discuss the case of $\PP(m)$: Proposition \ref{su2u1} says that 
we have a singular stabilizer for $m=3$ and a root stabilizer for $m=1$, while for all other values of $m$ the stabilizer
is generic regular; in all cases, except for $m=1$, we have that $N(U(1))$ is connected, hence we have $1$
possible way of gluing each of these tubular neighborood to others; for $m=1$ we have that the component $\tau\torus$
does not intersect $ S(U(2)\times U(1))$, so we have $2$ distinct classes in this case. The fact that the two classe 
obtained form the two gluing maps $\phi^{e}$ and $\phi^{\tau}$ cannot be isomorphic follows by inspecting diagram
(\ref{gluingdiag}). In fact all the vertical maps must be of the form $\phi^{e}$ in order to be defined on 
the whole tubular neighborhoods, and this implies that the central part of the diagram can not be commutative if we 
put, for instance, $\psi=\phi^{\tau}$.
In general we can combine two tubular neigborhoods if and only if the principal stanilizers are conjugate; therefore 
$\PP(n)$ and $\PP(m)$ can be glued together if and only if $n=m$, and the gluing map is unique for $m\neq1$, and there
are two distinct for $m=1$.\vspace{2mm}

 Let us pass now to tubular neighborhoods of type $\flag(p,q)$
assuming that $\gcd(p,q)=1$: there is precisely one gluing map for
$(p,q)\neq(0,1)$, $(p,q)\neq(1,0)$ or $(p,q)\neq(1,-1)$, as in fact in this case the
normalizers $N(U(1))$ are all connected. By contrast, for
the remaining representations the principal stabilizer is
of type $U_{1,-1}$, hence we have at first sight two possible gluing maps.
These can be used to join this tubular neighorhood to others with
the same type of stabilizer; nevertheless $\tau \torus\subset N(\torus)$, 
so as usual $\phi^{\tau}$ can be extended to the whole tubular neighborhood, and it is equivalent to the
identity gluing map. \vspace{2mm}

The list of all possible combinations is given in the Tables. $\blacksquare$
\vspace{4mm}

We end this section by examining the case in which $M/SU(3)$ is
$S^{1}$ and the principal orbits are simply connected. Let us point out that the homogeneous
manifold $\aw_{k,l}$ is simply connected, as shown by
the long exact homotopy sequence for a fibration:
\begin{equation}\label{lesaw}
\xymatrix{\cdots\,\pi_{1}(SU(3))\ar[r]&\pi_{1}(\aw_{k,l})\ar[r]&\pi_{0}(U(1))\,\cdots}\,.
\end{equation}
In this case there are no singular orbits and the manifold $M$ is a bundle
\begin{equation*}
\xymatrix{G/H\>\ar@{^{(}->}[r]&M\ar[d]\\ & S^{1}}
\end{equation*}
where $H=U_{k,l}$ is the principal (and unique) stabilizer; the structure group for this bundle is contained in $N(H)/H$ (see \cite[Th. 8.2, Ch. IV]{bredon}).
Hence we have

\begin{teor}\label{mgs1}
  Let $M$ be a cohomogeneity-one $SU(3)$-manifold with $M/G\cong
  S^{1}$ and such that the principal orbit is simply connected. Then
  the principal orbit has the form $\aw_{k,l}$.  Either
\begin{equation*}
M\cong \aw_{k,l}\times S^{1},
\end{equation*}
which is possible for any $k,l$, or $\aw_{k,l}=\aw_{1,-1}$ and $M$ is a nontrivial bundle over $S^{1}$.
\end{teor}

\proof We can divide the proof in three cases, corresponding to
the stabilizers described in Lemma \ref{normu1}.  First we note that
the bundle structure is given by the $N(H)/H$-valued transition
functions $g_{1}$ and $g_{2}$ defined on the two points $p_{1}$ and
$p_{2}$, which constitute the ``equator'' of the base manifold $S^{1}$.

In the first case, $N(U_{1,-1})/U_{1,-1}$ has $2$ connected components,
therefore there are two possible nonequivalent choices for the maps
$g_{i}$, giving rise to the trivial bundle and another nontrivial,
respectively.

In the remaining two cases, $N(U_{k,l})/U_{k,l}$ is connected: we have
a unique (trivial) bundle for $U_{1,1}$, and there 
are infinite nonconjugate generic $U_{k,l}$'s, giving rise to
nonisomorphic generic fibres $\aw_{k,l}$.

The $SU(3)$-manifolds obtained in this way are all trivial bundles,
except for the first case.  $\blacksquare$


\subsection{Non-connected singular stabilizers}\label{nonconnectedstabilizers}

We now conclude the classification, describing the more general
situation in which the singular stabilizers are not connected.  This
implies that the singular orbits are not simply connected, and their
respective universal covers are those described in Theorem
\ref{classu3}. Some of our arguments are inspired by those
used in \cite{alek-pod}.

\begin{prop}\label{propstab1}
If the connected components $K_{i}^{0}$ of the two singular stabilizers belong to the set $\{SO(3), SU(2), U(2)\}$, then
both are connected: $K_{i}^{0}=K_{i}$ for $k=1,\,2$.
\end{prop}

\proof\, Suppose that $K_{1}^{0}$ is one of the three subgroups in the
list: then the codimension of the singular orbit is at least $3$; a
general position argument shows that $M\setminus (SU(3)/K_{1})$ is
simply connected, as is $M$. This complement has the same homotopy
type of $SU(3)/K_{2}$, so $\pi_{1}(SU(3)/K_{2})=0$ too: this implies
that the stabilizer $K_{2}$ is connected. 
By the long exact homotopy sequence for a fibration
\begin{equation}\label{lesk1}
\xymatrix{\cdots\, \pi_{1}(S^{r})\ar[r]&\pi_{0}(H)\ar[r] &\pi_{0}(K_{2})\cdots}
\end{equation}
the principal stabilizer $H$ must also be connected, for $r>1$, which is the case for all
the representations involved with the three stabilizers under
consideration. This implies that also $K_{1}$ is connected, hence the result. $\blacksquare$ \vspace{2mm}

This means that we cannot obtain new simply-connected manifolds by
gluing together tubular neighborhoods unless they involve $\torus$ as
$K_{i}^{0}$ for at least one $i$. We discuss now this remaining case;
the new manifolds we obtain in this way are given in Table $4$. We
point out that the principal stabilizers turn out to be non-connected
in these cases. Before that, we prove a result which corresponds to
Lemma \ref{normu1} for non-connected $H$:

\begin{lemma}\label{normu1nc}
Consider the subgroup $U_{k,l}\times \Z_{h}$ of $\torus\subset SU(3)$; then 
$$N(U_{k,l}\times \Z_{h})=N(U_{k,l})$$
if $U_{k,l}$ is regular; if $U_{k,l}$ is singular (for instance $U_{1,1}$) then 
$$N( U_{1,1}\times \Z_{h})= \torus \cup \tau \torus\,.$$
\end{lemma}

\proof The proof in the regular case is completely analogous to that
of Lemma \ref{normu1}; for the singular case we just have to observe
that if $h\neq1$ the group contains regular elements, and the whole
$\torus$ must be preserved by the normalizer of $ U_{1,1} \times
\Z_{h}$. The element $\tau\in\weyl$ that normalizes $U_{1,-1}$,
reflecting the root $\alpha$, is the only one which also preserves
$U_{1,1}\times \Z_{h}$, hence the conclusion.
$\blacksquare$\vspace{2mm}

\noindent\emph{Observation.} In this situation, we have two connected components for the normalizers of 
$U_{1,-1}\times \Z_{h}$ and of $U_{1,-1}\times\Z_{h}$; if these stabilizers appear in a tubular
neighborhood of type $\flag(p,q)$, we observe that in both cases we
obtain only one $SU(3)$ diffeomorphism class, because both normalizers
are contained in $N(\torus)$ (for Uchida's criteria, see Theorem
\ref{classu3}).\vspace{2mm}

We pass now to the main result of this section, but before of that we recall that any subgroup $K\subset G$ is 
alaways contained in $N(K^{0})$, because for any $x\in K$ the adjoint action $Ad_{x}$ is continuous, preserves $K$ 
and fixes $e$.

\begin{prop}\label{propstab2}
Suppose that $K_{1}^{0}\in\{SO(3), SU(2), U(2)\}$ and that $K_{2}^{0}=\torus$: then $K_{2}=K_{2}^{0}=\torus$.
Moreover if $K_{2}=\torus$ and if the slice representation at $\flag$ is $A^{p}\otimes A^{q}$ there are the following possibilities:
\begin{enumerate}
\item
if $(p,q)=(0,h)$ then $K_{1}\in\{SO(3),\,SU(2),\,U(2),\, \torus\}$;
\item
if $(p,q)=(0,h)$ for some $h\in\Z,\,h>1$, then $K_{1}=(SU(2)\times
\Z_{2h})/\Z_{2}$, except in the case $h=3$, where 
also $K_{1}=SO(3)\times \Z_{3}$ is possible;
\item
if $(p,q)\neq(0,h)$ and $\gcd(p,q)=1$ then $K_{1}\in\{\torus, \,U(2)\}$;
\item
if $(p,q)\neq(0,h)$ and $\gcd(p,q)\neq1$ then $K_{1}^{0}=\torus$.
\end{enumerate}
\end{prop}

\begin{table}[h]
\begin{center}
\setlength{\extrarowheight}{1mm}
\setlength{\tabcolsep}{2mm}
\vspace{0.5cm}
\begin{tabular}{|c||c|c|c|}
\hline
 $_{M_{2} }\backslash^{M_{1}}$  & $\SS/\Z_{h}$ & $\ssu/\Z_{3}$ & $\flag(l,m)$ \\[1pt]
\hline
\hline
  $\flag(p,q)$ & $\delta^{p}_{0}\delta^{h}_{q}$ & $\delta^{p}_{0}\delta^{3}_{q}$ & $\delta^{p}_{l}\delta^{m}_{q}$ \\[1pt]
\hline
\end{tabular}
\parbox{350pt}{\caption{Numbers of $SU(3)$-diffeomorphism classes of
    $8$-manifolds: non-connected principal stabilizers ($\gcd(p,q)\neq1$)}}
\end{center}
\label{nonconn}
\end{table}

\proof The first statement follows from the same general position
argument as in Proposition \ref{propstab1}. As proved in Lemma
\ref{lemmastabtorus}, the principal stabilizer is of the form
$U_{k,l}\times \Z_{h}$, so we need to determine which of the
singular stabilizers contain this subgroup.

In the first case we have $h=1$ and $U_{k,l}=U_{1,-1}$, which appears as a
principal stabilizer associated to any of the connected stabilizers
above, with the appropriate slice representation $\slice$, as already
shown in Theorem \ref{classu3}. In this case $K_{1}$ is connected,
because the sphere $S^{r}\subset\slice$ is.

For the second case we argue as follows: $K_{1}$ must contain a subgroup of type 
$U_{1,-1}\times \Z_{h}$, but we have to exclude $U(2)$, because it allows
only connected principal stbilizers ($h=1$). 
Another possible choice is the subgroup 
\begin{equation*}
\frac{SU(2)\times \Z_{2h}}{\Z_{2}}
\end{equation*}
where $\Z_{2}$ is the center of $SU(2)$; topologically it is the union of $h$ copies of
$SU(2)$, and $\Z_{2h}$ should be regarded as a subgroup of the singular $U(1)$
centralizing $SU(2)$ (for instance $U_{1,1}$ for the standard immersion).
We observe that the singular orbit in this case is isomorphic to $\SS/\Z_{h}$. 
Suppose instead that $K_{1}^{0}=SO(3)$: then $K_{1}$ must be a subgroup of
the normalizer $N(SO(3))=SO(3)\times \Z_{3}$, which are $SO(3)$ itself or 
the whole $N(SO(3))$. The latter case in this situation corresponds for instance to the
weight $(p,q)=(0,3)$ for $\torus$. As observed after Lemma \ref{normu1nc},
in both cases the two gluing maps $\phi^{e},\psi^{\tau}$ give rise to isomorphic 
$SU(3)$-spaces.

For the third case, the connected component $K_{1}^{0}$ must contain
$\torus$, because only then the corresponding Lie algebra does contain
the correct $\uuu(1)$. 

Finally, in the fourth case we have to exclude $U(2)$ because, as observed in case $2$, it allows only connected 
principal stabilizers. $\blacksquare$
\vspace{2mm}

\noindent\emph{Observation.} Two of the manifolds that are new with respect to the classification given in 
Theorem \ref{classu3} come from case $2$. 
We note that in these cases the singular stabilizer $(SU(2)\times \Z_{2h})/\Z_{2}$ admits as a slice representation $\slice$
only the standard $\R^{3}\cong [\si^{2}\otimes A^{0}]$: in fact any $\Z_{h}$ representation can be extended 
to a $U(1)$ representation $A^{m}$ with $0\leq m\leq h-1$, hence $\slice$ is the restriction of a $U(2)$ representation; 
therefore 
\begin{equation}
\slice_{\C}\cong \sum \si^{l}\otimes A^{m}
\end{equation}
as seen in Proposition \ref{su2u1}; for dimensional reasons $\si^{2}\otimes A^{0}$ is the 
only possible choice. Analogous considerations hold for $SO(3)\times \Z_{3}$. 
\vspace{2mm}

Let us consider the case in which $K_{i}\neq K_{i}^{0}=\torus$, so that $K_{i}\subset N(\torus)$. Here,
the two stabilizers must have the same number of connected components,
otherwise the two tubular neighborhoods could not be glued together, as the principal orbits would 
not be isomorphic. In this case $K_{1}=K_{2}$ and 
$\pi_{1}(SU(3)/K_{i})\neq0$ for $i=1,2$; moreover in the long exact sequence 
\begin{equation*}
\xymatrix{\cdots\, \pi_{1}(SU(3)/H)\ar[r] &\pi_{1}(SU(3)/K_{i})\ar[r]&\pi_{0}(S^{1}) \,\cdots}
\end{equation*}
the bundle pro\-jec\-tions induce sur\-jec\-tions on the res\-pec\-tive
fun\-da\-men\-tal groups. The Seifert--van Kampen Theorem tells us that this
is incompatible with the simply connectedness of the manifold $M$,
so we have to exclude this case.  \vspace{2mm}


\section{Examples}\label{examplesandapplications}

In order to present some familar examples, the notation 
$\sM(M_{1},\,M_{2})$ indicates an $8$-dimensional $SU(3)$-manifold
obtained by gluing appropriate disk bundles over singular orbits
$M_{1},\,M_{2}$ with a map $\phi$ which may or may not be the
identity. 

Then we have the following remarkable identifications:
\begin{itemize}
\item
the complex Grassmannian $\grasc{4}$ is $\sM(\PP,\,\PP)$;
\item
the quaternionic projective plane $\Ha\Pro^{2}$ is $\sM(\PP,\,\SS)$
\item
the exceptional Wolf space $G_{2}/SO(4)$ is $\sM(\PP,\,\ssu)$
\item the product $\C\Pro^2\times\C\Pro^2=\PP\times\PP$ is
  $\sM(\PP,\flag)$
\item
the Lie group $SU(3)$ is itself $\sM(\ssu,\,\SS)$.
\end{itemize}

We describe these $SU(3)$ spaces in a bit more detail. Recall that
$\ssu=SU(3)/SO(3)$.  The first three examples are obtained by standard
inclusions of $SU(3)$ in $SU(4)$, $Sp(3)$ and $G_{2}$, and in these
cases, the normal bundle over each $\C\Pro^2=\PP$ is $\PP(1)$.

The fourth (product) case is given by the diagonal action of
$SU(3)$, where the first singular orbit consists of all couples
$([z],[z])$ of identical complex lines in $\C^{3}$, and the second
consists of couples $([z],[w])$ with $[w]\subset [z]^{\perp}$.  In
this case, the slice representation $\slice$ is isomorphic to the
isotropy representation at $\PP$: in fact if $(v,v)$ is a tangent
vector at $([z],[z])$, with $v$ generated by an elelement in
$\sut/(\uuu(2)\oplus\R)$, then normal vectors must be of the form
$(v,-v)$ and give rise to the same $U(2)$ representation. It is
straightforward to check that this tubular neighborhood is of type
$\PP(3)$.

The final case is given by a modification of the Adjoint action of
$SU(3)$ on itself, discussed in more detail in
Subsection~\ref{theconsimilarityaction}.\vspace{2mm}

The case in which the two tubular neighborhoods are isomorphic and the
gluing map is the identity is particularly simple. We can identify the
singular orbits $M_{1}=M_{2}=M$ and call the unique normal
representation $\slice$; the resulting manifold
$\double(M)=\sM(M_1,M_2)$ is then the ``double'' of the disk bundle
associated to $V$. This manifold is obtained by the one-point
compactification $\R^{n}\leadsto S^{n}$ of the $\slice$ fibres over
$M$:
\begin{equation}\label{dou}
\xymatrix{S^{n}\>\ar@{^{(}->}[r]& \double(M)\ar[d] \\ & M.}
\end{equation}
The other singular orbit becomes the section at infinity of this new
bundle. \vspace{1mm}

\begin{prop}
The manifolds $\double(\SS)$ and $\double( {\ssu})$ do not admit any $SU(3)$-invariant metric of positive
sectional curvature.
\end{prop}
\proof This is just a consequence of \cite[Lemma 3.2]{podestaverdiani2}, which asserts that any even dimensional 
cohomogeneity-one $G$-manifold $M$ with an invariant metric of positive sectional curvature has $\chi(M)>0$, and
of the observation 
\begin{equation*}
\chi(\double(\SS))=\chi(\double(\ssu))=\chi(S^{3}) \chi(\SS)=0
\end{equation*}
(recall that $\chi(\ssu)=\chi(\SS)$). $\blacksquare$


\subsection{Consimilarity}\label{theconsimilarityaction}

We are going now to consider a group action $\mathbf{c}$ of
$GL(n,\C)$ on itself, called \emph{consimilarity}, defined by
\begin{equation}\label{consim}
  \cons(A)B := A\kern1pt B\kern1pt \overline{A}^{-1}.
\end{equation}
This action naturally occurs when considering \emph{anti}-linear
mappings between a given vector space, of relevance in quantum theory.
It also occurs in various geometrical situations (see, for example,
\cite{finopartsala}), and is intimately related to \emph{similarity}.
The mapping
\begin{equation}\label{AA}
  \xxbar\colon A\mapsto A\kern1pt \overline A
\end{equation} induces a mapping between 
consimilarity classes and similarity classes (i.e.\ 
orbits under (\ref{consim}) and orbits under conjugation).
Although this mapping is not in general a bijection between the
respective classes, it is true that $\xxbar^{-1}(I)$
coincides with the consimilarity orbit 
\begin{equation}\label{AAminus}
\{A\kern1pt \overline A^{-1}:A\in GL(n,\C)\}
\end{equation}
of the identity. This fact is not entirely obvious, but has an easy
proof \cite{hornjohnson}. 

Consimilarity can be restricted to $SU(n)\subset GL(n,\C)$,\,so that
$SU(n)$ acts on itself, as in this case
\begin{equation*}
A\kern1pt B\kern1pt \overline{A}^{-1}=AB A^{t}
\end{equation*}
is in $SU(n)$ if $A,\,B$ are.  It is straightforward to prove that the
consimilarity action of $SU(n)$ on itself is isometric with respect to
the Killing metric. The resulting action is in fact a special case
of a family of actions of a Lie group $G$ on itself, constructed using
an automorphism $\sigma$ of $G$ (see \cite{hptt1}, \cite{conlon} and
\cite{hornjohnson}).

Let us return to the case $n=3$. 

\begin{lemma}\label{consimlemma}
  Consimilarity is a cohomogeneity-one action of $SU(3)$ on itself
  with singular orbits $\ssu$ and $\SS=S^5$. The former is the orbit
  containing the identity matrix $I$ and coincides with
  $\xxbar^{-1}(I)\cap SU(3)$.
\end{lemma}

\proof It can be shown that $\xxbar^{-1}(I)\cap SU(3)$
coincides with the set \[\mathcal{S}=\{A\in SU(3):A=A^t\}\] of
symmetric matrices.  The map $\xi:\ssu\to\mathcal{S}$ defined by
\[A \,SO(3)\xymatrix{\ar@{|->}[r] &} AA^{t}\] is well defined and
surjective. It is also injective as if $AA^{t}=CC^{t}$ then
\begin{equation*}
C^{-1}A=C^{t}(A^{t})^{-1}=\big((C^{-1}A)^{t}\big)^{-1}
\end{equation*}
so that $C^{-1}A\in SO(3)$. This shows that the $\cons$-orbit through
the identity is $\ssu$.  

Consider the point $I$ and its stabilizer $SO(3)$. The isotropy and
the slice representations are determined by the decomposition
\begin{equation*}
\sut=T\ssu\oplus\slice=\sot^{\perp}\oplus\sot=[\si^4]\oplus [\si^2]
\end{equation*}
as $SO(3)$ representations. The slice representation is sphere transitive (see Table~1),  
and this shows that the cohomogeneity of the action is $1$.
For instance we can choose the normal direction determined by the matrix
\begin{equation*}
w=
\begin{pmatrix}
0&1&0\\
-1&0&0\\
0&0&0
\end{pmatrix}
\in \sot=\si^2\,.
\end{equation*} 
The corresponding geodesic $B(t)=\exp(tw)$ intersects orthogonally all
the $\cons$-orbits (see \cite{hptt1}): the second singular one is
reached at $B_{s}=B(\pi/4)$. In fact, an explicit calculation shows
that the stabilizer at $B_{s}$ is $SU(2)$; therefore the corresponding
orbit is indeed $\SS$. $\blacksquare$ \vspace{2mm}

Observe that
\begin{equation}\label{imtrz}
\overline{\Tr(A\kern1pt \overline{A})}=\Tr(\overline A\kern1pt{A})
=\Tr(A\kern1pt \overline{A})\;;
\end{equation}
this implies that the image $\xxbar(SU(3))$ is contained in the
hypersurface 
\begin{equation}\label{hyper}
\hypers:=\{B\in SU(3):\Tr\,B\in\R\}
\end{equation}
of $SU(3)$. We shall investigate the resulting mapping $\xxbar\colon
SU(3)\to \hypers$ in the next section.


\section{Quotients by circle subgroups}\label{quotientsbycirclesubgroups}

An analogous classification of $SU(3)$ actions is possible in
dimension $7$, and partial results can be found in
\cite{podestaverdiani}. Restricting attention here to the case in
which both singular orbits are $\PP=\C\Pro^2$, and both tubular
neighborhoods are isomorphic to the rank 3 vector bundle
$\Lambda^{2}_{-}\C\Pro^2$, it is not hard to show the existence of only two
classes of cohomogeneity-one $SU(3)$-spaces with this data. There is a
choice of gluing map between the generic orbits $\flag$: the identity
in one case, and a map $\phi^{\tau}$ associated to a non-trivial
element $\tau\in\weyl$ in the other. With the latter choice, we
obtain the sphere $S^{7}\subset\sut$ with the action induced by the
Adjoint representation.

We now exhibit a model for the manifold obtained in the former case,
denoted here by $N^{7}$, involving the Grassmannian $\graor(\sut)$ of
oriented 3-dimensional subspaces of the Lie algebra $\sut$, which is 
an $SU(3)$-space under the action induced by $Ad_{SU(3)}$.

\begin{prop}\label{N7}
  The manifold $N^{7}$ is a submanifold of $\graor(\sut)$
  with the $SU(3)$ action induced by the Adjoint action on $\sut$.
\end{prop}

\proof Following \cite{swann98}, we consider the function $f\colon
\graor(\sut)\to\R$ induced by the standard $3$-form on $\sut$.  Thus
$$f(U)=\langle x,[y,z]\rangle,$$
where $\{x,y,z\}$ is an orthonormal
basis of the $3$-dimensional subspace $U\subset\sut$.  The absolute
maxima and minima of $f$ are each attained on a copy of $\C\Pro^2$
corresponding to the highest root embedding $\sud\subset\sut$ and a
choice of orientation for $\sud$. The tangent space $T_U\graor(\sut)$
has the form $U\otimes U^{\perp}$; for $U=\sud$ it can be decomposed
as
\begin{equation*}
T_{\sud}\graor(\sut)=\sud\otimes (\si^0 +2\si^1)\cong \si^2+ 2(\si^1+\si^3).
\end{equation*}  
The subspace $2\si^{1}=\si^{1}\oplus\si^1$ represents the tangent space to
the critical manifold $\C\Pro^2$; if we choose instead the summand
$\si^2$, we obtain the bundle $\Lambda^{2}_{-}\C\Pro^2$, which is therefore a subbundle 
of the normal bundle at both $\C\Pro^{2}$. In the two cases it turns out to be a stable 
or an unstable subbundle respectively. 

The manifold $N^7$ is obtained from the two $\Lambda^{2}_{-}\C\Pro^2$ over the
two extremal $\C\Pro^2$. To see this, denote by $\tilde{N}^{7}$ the manifold
obtained by considering the union of the flow lines of the vector field $\grad\,f$  
with limit points in the two copies of $\C\Pro^{2}$ and tangent directions corresponding to the 
respective $\si^{2}$. Such a flow line (without caring about the parametrization) is given by
\begin{equation*}
V(t)=\vspan\{\ugen\cos{t}+\vgen\sin{t},\,\ugen_{2},\,\ugen_{3}\}\,,
\end{equation*}
with $\ugen,\vgen$ as in (\ref{princu1}) and $\sud=\vspan\{\ugen
,\ugen_{2},\ugen_{3}\}$.  It is straightforward to see that the 
stabilizer for $t\neq k\pi$ under the $Ad_{SU(3)}$ action is $\torus$, and for
$t=\pi$ the integral curve intersects the minimal critical
submanifold at the same subalgebra $\sud$ with opposite orientation.
In both cases the tangential direction of $V(t)$ belongs to the
summand $\si^2$ at the critical points: these facts imply that the gluing map for the two
tubular neighborhoods must be the identity, so
$\tilde{N}^{7}\cong N^{7}$. $\blacksquare$ \vspace{3mm}

\noindent \emph{Remark.} We point out that $N^{7}$ is \emph{not}
homeomorphic to $S^{7}$. As the double $\double(\C\Pro^2)$, it can
be regarded as a $3$-sphere bundle over $\C\Pro^2$ as in (\ref{dou}), in
contrast to $S^7$. Now, $\pi_{2}(\C\Pro^2)=H_{2}(\C\Pro^2,\Z)=\Z$, and writing
the homotopy exact sequence for a fibration we obtain
\begin{equation*}
\xymatrix{\cdots\, \pi_{2}(S^{3})\ar[r] &\pi_{2}(N^{7})\ar[r]&\pi_{2}(\C\Pro^2)\ar[r]&\pi_{1}(S^{3}) \,\cdots}.
\end{equation*}
This implies $\pi_{2}(N^{7})=\pi_{2}(\C\Pro^2)=\Z$, whilst $\pi_{2}(S^{7})=0$.

It is shown in \cite{podestaverdiani} that $N^{7}$
cannot be equipped with an invariant metric of positive curvature.
Indeed, $S^{7}$ is the unique $7$-dimensional positively curved
cohomogeneity-one $G$-manifold, if the semisimple part of $G$ has
dimension greater then $6$. 
\vspace{3mm}

The above example is linked to the $8$-dimensional case by a moment
map $\mu$ associated to the action of $SU(3)$ on the 
Wolf spaces $\Ha\Pro^{2}$ and $\grasc{4}$ (see \cite{galaw1}). Denoting by $M$ either of
these space, it is possible to construct from $\mu$ an equivariant map 
\begin{equation}\label{mmapsu3}
\Psi: M_{0} \xymatrix{\ar[r]&}\graor(\sut)
\end{equation}
defined on an open dense subset $M_{0}\subset M$. This construction
was used in \cite{gambioli1} in order to relate the geometry of a
quaternion-K\"{a}hler manifold with the geometry of the Grassmannian
$\graor(\la)$, but we cannot use the same techniques here since in the
two cases considered, the differential $\Psi_{*}$ is nowhere injective.
Moreover the subset $M_{0}$ is strictly contained in $M$:
indeed $M_{0}= \Ha\Pro^{2}\setminus \C\Pro^2$ and $M_{0}=\grasc{4}\setminus
\C\Pro^2$ respectively, and
\begin{equation*}
\Psi(M_{0})\subset N^{7}\subset \graor(\sut)\,.
\end{equation*}

One may ask if the map $\Psi$ could be extended equivariantly to the
whole $W$ in both cases, as this happens in other significant cases
(for instance $Sp(n)Sp(1)$ acting on $\Ha\Pro^{n}$ or $Sp(n)$ acting
on $\grasc{n}$).  In fact the generic fibre $\Psi^{-1}(x)$ is a circle
$S^{1}$: the resulting $S^{1}$ action on $\Ha\Pro^{2}$ was described
in \cite{batt}, and 
\begin{equation*}
  \Ha\Pro^{2}/S^{1}\cong S^{7}
\end{equation*}
(see \cite{atiyahberndt} and \cite{atiyahwitten}). For the same reason we have a
topological quotient
\begin{equation*}
\grasc{4}/S^{1}\cong S^{7}.
\end{equation*}
However, as observed above, $S^{7}$ is different from $N^{7}$, and it
is easy to check that $\Psi$ cannot be extended equivariantly to the
whole Wolf spaces $\Ha\Pro^{2}$ and $\grasc{4}$.

The fact that $S^7$ is a compactification of $\Lambda^2_{-}\C\Pro^2$ was
used in \cite{miyaoka} to unify the construction of various Ricci-flat
metrics on complements of homogeneous spaces inside spheres. There
are analogous constructions on $G_2/SO(4)$ and $SU(3)$. The
descriptions at the start of Section \ref{examplesandapplications} show that dense open subsets
of thse two manifolds can be $SU(3)$-equivariantly identified.
However, the respective singular orbits $\PP$ and $\SS$ are not
directly related by the Hopf fibration $\SS\to\PP$; indeed passing
from the $\PP$ of $G_2/SO(4)$ to the $\SS$ of $SU(3)$ requires a
``flip'' of the type considered in \cite{gukospartong}. This is made
possible by the existence of three distinct mappings $\flag\to\PP$,
similarly exploited in the theory of harmonic maps \cite{sal85}.
\vspace{1mm}

To conclude the paper, we identify an analogue for $SU(3)$ of the
map $\Psi$ described in (\ref{mmapsu3}). 

\begin{teor} The image $\xxbar(SU(3))$ of (\ref{AA}) is the
  hypersurface (\ref{hyper}), and is homeomorphic to the Thom space
  of the vector bundle $\Lambda^2_-\C\Pro^2$. The restriction of $\xxbar$
  to $SU(3)\setminus\ssu$ is a principal $S^1$ bundle over
  $\Lambda^2_-\C\Pro^2$.
\end{teor}

\proof

Let $SU(2)\subset SU(3)$; then
\begin{equation}\label{eqimtr}
\hypers=\bigcup_{g\in SU(3)} Ad_{g}SU(2).
\end{equation}
In fact, consider the Lie algebra $\sut$; it is well known
that any element $x\in \sut$ belongs to the standard $\ta=\vspan\{\ugen,\vgen\}$ (see (\ref{basiscartan})), up
to conjugation. It is sufficient therefore to solve the equation
\begin{equation*}
\mathrm{Im}\,\Tr\big(\exp{(t\,\ugen+s\,\vgen)}\big)=0,
\end{equation*}
equivalent to
\begin{equation*}
\sin(t+s)+\sin(s-t)-\sin(2s)=0\,;
\end{equation*}
this has solutions $\{s=0+k\pi\} \cup \{s=\pm t +2k\pi\}$. These are
nothing other than the three lines corresponding to $\uo{1,-1}$,
$\uo{1,0}$, $\uo{0,1}$ and their translates. However, when we
exponentiate, all the solutions are sent to the triplet 
\begin{equation}\label{toruscaph}
U_{1,-1}\cup U_{1,0}\cup U_{0,1}=\torus \cap \hypers 
\end{equation}
which are the intersections of $\torus$ with
three conjugate copies of $SU(2)$; the equality in (\ref{eqimtr})
follows by noting that both sides are $Ad$-invariant.
Consider again any subgroup $SU(2)\subset SU(3)$ and let $g\in
SU(3)$: then
\begin{equation}\label{sudint}
SU(2) \cap Ad_{g}SU(2)= SU(2)\hbox{ or }\{e\}.
\end{equation}
In fact let us consider a point $x\in SU(2)\cap Ad_{g}SU(2)$; if $x$ is regular then it is contained in a unique 
maximal torus $\torus$; on the other hand $x$ belongs to one of the connected components of (\ref{toruscaph}), say $U_{1,-1}$, which therefore
belongs entirely to $SU(2)\cap Ad_{g}SU(2)$. This implies that $g\in N(U_{1,-1})$, which is contained 
in $N(SU(2))$, hence we fall in the first case of (\ref{sudint}). Suppose now that $x$ is singular: there exists only
one singular point for each copy of $SU(2)$, namely 
\begin{equation}\label{antipsud}
x=
\begin{pmatrix}
-1&0&0\\
0&-1&0\\
0&0&1
\end{pmatrix}\,
\end{equation}
for the standard embedding of $SU(2)$; singular elements are preserved
by the adjoint action, therefore $g$ is in the stabilizer of $x$,
which is $U(2)=N(SU(2))$, and we are again in in the first case of
(\ref{sudint}). If $g\not\in N(SU(2))$ then the intersection consists
of just $e$.

This discussion proves that we can realize $\hypers$ as the union of
copies of $SU(2)$ which share only the identity $e$ inside $SU(3)$. On
the other hand, the singular orbit $\C\Pro^2$ parametrizes this union; our
conclusion is that $\hypers$ is therefore isomorphic to the total space of a
fibre bundle $\princ$ over $\C\Pro^2$ with fibre $SU(2)$ and with one point
for each fibre identified:
\begin{equation}
\xymatrix{ \princ\ar[d]& \> SU(2)\ar@{_{(}->}[l] \\
  \C\Pro^{2}&}
\end{equation}
and $\hypers=\princ /\sim$, with $e\sim e'$ if and only if $e$ and $e'$ are the identity of two fibres $SU(2)$ and $SU(2)'$ (the identity is well defined as it 
is fixed by the isotropy subgroup of $\C\Pro^2$ acting on the fibres). 

The Thom space of a vector bundle $E\to M$ is obtained by a $1$-point
compactification of the total space $E$. Our construction shows that
$\hypers$ is indeed the Thom space of the bundle
$\Lambda^{2}_{-}\C\Pro^2$: in fact the fibre of this vector bundle is
isomorphic to $\sud$ as a representation of the stabilizer $U(2)$;
then, consider the closed disk $D_{\!\sqrt{2}\pi}\subset \sud$: we can
identify
\begin{equation*}
SU(3)=\exp D_{\!\sqrt{2}\pi} 
\end{equation*}
where the spheres $S^{2}_{r}$ of radius $r<\sqrt{2}\pi$ are sent 
$Ad_{SU(2)}$-equivariantly to $2$-spheres, whilst the boundary $S^{2}_{\sqrt{2}\pi}$
is collapsed to a point $x$ antipodal to $e$ (see (\ref{antipsud})). We have therefore
a corresponding disk subbndle $D$, and a bundle with fibre $SU(2)\cong S^{3}$
obtained from the former by collapsing the boundary of each fibre to a point.
The Thom space can be therefore obtained by additionally identifying all
the antipodal points of the various fibres. This is precisely what happens for the
hypersurface $\hypers$, but this time identifying the identities $e$ instead of the antipodal
points. This is not a real difference: in fact the antipodal element $x$ belongs to the center 
$C(SU(2))=\Z_{2}$, and the automorphism $SU(2)\to xSU(2)$ 
is $Ad_{SU(2)}$ equivariant and swaps $e$ and $x$, giving rise to isomorphic bundles with
fibre $SU(2)$. The hypersurface $\hypers$ can be shown to be smooth
everywhere excepted at $e$.

The image $\xxbar(SU(3))$ is contained in $\hypers$, as seen in
(\ref{imtrz}); the surjectivity of $\xxbar$ can be established in the
following way by equivariance: the normal geodesic $B(t)$ used in
Lemma \ref{consimlemma} intersects all the $\cons$ orbits; its image
is given by
\begin{equation*}
\xxbar(B(t))=B(2t)
\end{equation*}
which intersects all the $Ad_{SU(3)}$ orbits orthogonally, joining the two singular orbits $e$ and $\C\Pro^2$. 
We observe that the singular orbit $\ssu\subset SU(3)$ is collapsed to $e$.

We pass now to the last statement of the theorem: we will use an
argument which is a bundle version of that discussed in Section
\ref{preliminaryresults} (see (\ref{equivhom})). We can describe the
tubular neighborhood $D_{\SS}\cong SU(3) \setminus \ssu$ around $\SS$
as the $[\si^{2}]=\R^{3}$ bundle obtained by the twisted product
\begin{equation*}
SU(3)\times _{\scriptscriptstyle SU(2)} [\si^{2}]\,;
\end{equation*}
in other words the couples $(g,v)\in SU(3)\times \R^{3}$ are
identified by the relation
$(g,v)\sim (g',v')$ if and only if $ g'=g h\;,\,v= h^{-1}v'$
for some $h\in SU(2)$. The space of classes $[g,v]$ is naturally a left $SU(3)$-space under the action
$g'[g,v]=[g'g,v]$. Observe now that the $SU(2)$ representation $\R^{3}$ can be extended to a $U(2)$ 
representation of the form
$[\si^{2}]\otimes A^{0}$, so that the $U(1)$ centralizer of $SU(2)$ acts trivially. This implies that  
$D_{\SS}$ becomes also a \emph{right} $U(2)$-space in the following way: an element $k\in U(2)$ acts by
$k[g,v]=[gk,k^{-1}v]$. This action is well defined because $U(2)=N (SU(2))$, and it is equivariant with respect to the left $SU(3)$ action.
Clearly $SU(2)\subset U(2)$ is precisely the non-effectivity kernel, 
so we can just consider this action a $U(2)/SU(2)=U(1)$ effective action. The quotient space $D_{\SS}/U(1)$ turns out
to be a twisted product of the form $SU(3)\times _{\scriptscriptstyle U(2)} \slice$,
with $\slice=[\si^{2}]\otimes A^{0}$, which is nothing other than
$\Lambda^{2}_{-}\C\Pro^2$. The projection $\pi_{U(1)}$ is therefore an
equivariant map
\begin{equation}
\xymatrix{SU(3)\setminus \ssu\ar[r]& \hypers\setminus\{e\}}
\end{equation}
as is the map $\xxbar$; the restriction to each orbit is an
equivariant projection of homogeneous spaces in both cases, and an
inspection of the normalizers of $SU(2)$ and $U_{1,-1}$ shows that
the choice is unique, hence $\xxbar=\pi_{ U(1)}$. $\blacksquare$
\vspace{2mm}

\noindent\emph{Observation.} The proof above has identified $\xxbar$
with the quotient 
\begin{equation*}
SU(3)\setminus\ssu\>\cong\>
\xymatrix{\Ha\Pro^{2}\setminus \C\Pro^2\ar[r]& 
S^{7}\setminus \C\Pro^{2}\cong \Lambda ^{2}_{-}\C\Pro^{2}}\,.
\end{equation*}
induced by the $U(1)$ action described in \cite{atiyahwitten}, \cite{miyaoka}.
\vspace{2mm}

Complete metrics of holonomy $Spin(7)$, invariant under a $Spin(5)$
action, have been discovered on the positive spin bundle over $S^{4}$
\cite{brysal}; more recently other metrics of this type have been
constructed on $4$-di\-men\-sio\-nal vector bundles over $\C\Pro^{2}$
(see \cite{gukovsparks}). These bundles belong to the family we have
denoted by $\PP(l)$ (see Proposition \ref{su2u1}). In a future
article, we hope to use the examples of this paper to construct new
special geometries in dimensions 7 and 8, by gluing together tubular
neighborhoods that arise in our classification, adapting invariant
structures to appropriate conditions at the boundaries.\vspace{5mm}

\noindent{\footnotesize{\bfseries Acknowledgements.} The author wishes to thank
  S. Salamon for his constant support and encouragment. He is also
  grateful to Y. Nagatomo and F. Podest\`a for essential help in
  getting this paper underway, and to F. Lonegro, M. Pontecorvo and A.
  Di Scala for additional input. The paper was written whilst the
  author was a recipient of a grant within the research projects ``Geometria 
 Riemanniana e strutture differenziabili" (University of Rome \emph{La Sapienza}) 
and ``Geometria delle variet\ac differenziabili" (University of Florence).}

\vspace{4mm}

\noindent\textsc{Dipartimento di Matematica ``G. Castelnuovo", 
Universit\ac ``La Sapienza", Piazzale A. Moro, 2 - 00185 Roma - Italy }\\
\emph{E-mail address:} \texttt{gambioli@mat.uniroma1.it}


\begin{thebibliography}{99}\frenchspacing

\bibitem{alek-alek} A.V. Alekseevsky, D.V Alekseevsky: Riemannian
  $G$-manifold with one-dimensional orbit space, Ann. Global Anal.
  Geom. \textbf{11} (1993), 197--211.

\bibitem{alek-alek2} A.V. Alekseevsky, D.V Alekseevsky: $G$-manifold
  with one-dimensional orbit space, Adv. in Sov. Mat. \textbf{8} (1992), 1--31.

\bibitem{alek3} D.V. Alekseevsky: Compact quaternion spaces, 
Functional Anal. Appl., \textbf{2} (1968), 106--114.

\bibitem{alek-pod} D. V. Alekseevsky, F. Podest\`a: Compact
  cohomogeneity one Riemannian manifolds of positive Euler
  characteristic and quaternionic K\"{a}hler manifolds, Geometry,
  Topology, Physics. Proceedings of the First USA-Brazil Workshop,
  Campinas 1996 (B. N. Apanasov et al. eds.), de Gruyter, Berlin,
  1997, 1--33.

\bibitem{aw} S. Aloff, N. Wallach: An infinite family of distinct
  $7$-manifolds admitting positively curved Riemannian structures,
  Bull. A.M.S. \textbf{81} (1975), 93--97.

\bibitem{atiyahberndt} M. Atiyah, J. Berndt: Projective planes, Severi
  varieties and spheres, Surveys in Differential Geometry VIII, Papers
  in Honor of Calabi, Lawson, Siu and Uhlenbeck, International Press,
  Somerville, MA, 2003, 1--27.

\bibitem{atiyahwitten} M. Atiyah, E. Witten: M-theory dynamics on a
  manifold of $G_{2}$ holonomy, Adv. Theor. Math. Phys. \textbf{6}
  (2002) 1--106

\bibitem{batt} F. Battaglia: Circle actions and Morse theory on
  quaternion-K\"ahler manifolds, J. London Math. Soc. \textbf{59}
  (1999), 345--358.

\bibitem{beber} L. B\'{e}rard Bergery: Sur de nouvelles
  vari\'{e}t\'{e}s riemanniennes d'Einstein, Publications de
  l'Institut E. Cartan \textbf{4} (Nancy, 1982), 1--60.

\bibitem{besse} A. Besse: \emph{Einstein Manifolds}, Springer-Verlag,
  1987.

\bibitem{bredon} G.E. Bredon: \emph{Introduction to compact
    transformation groups}, Number \textbf{46} in Pure and Applied
  Mathematics, Academic Press, 1972.

\bibitem{brocktd} T. Br\"{o}cker, T. tom Dieck: \emph{Representations
    of Compact Lie Groups}, Sprin\-ger, 1985.

\bibitem{brysal} R.L. Bryant, S.M. Salamon: On the construction of
  some complete metrics with exceptional holonomy, Duke Math. J.
  \textbf{58} (1989), 829--850.

\bibitem{conlon} L. Conlon:  The topology of certain spaces of paths
  on a compact symmetric space, Trans. Amer. Math. Soc. \textbf{112} (1964), 228--248.

\bibitem{dancswan1} A. Dancer, A.F. Swann: Quaternionic K\"{a}hler
  manifolds of cohomogeneity one, Int. J. Math. \textbf{10} (1999),
  541--570.

\bibitem{dancswan2} A. Dancer, A.F. Swann: Hyperk\"{a}hler metrics of
  cohomogeneity one, J. Geom. and Phys. \textbf{21} (1997), 218--230.

\bibitem{dancwang1} A. Dancer, M.Y. Wang: Painlev\'e expansions,
  cohomogeneity one metrics and exceptional holonomy.  Comm. Anal.
  Geom. \textbf{12} (2004), 887--926.

\bibitem{finopartsala} A. Fino, M. Parton, S. Salamon: Families of
  strong KT structures in six dimensions, Comment. Math. Helv.
  \textbf{79} (2004), 317--340.

\bibitem{galaw1} K. Galicki, B. Lawson: Quaternionic reduction and
  quaternionic orbifolds, Mat. Ann. \textbf{282}(1988), 1--21.
  
\bibitem{gambioli1} A. Gambioli: Latent quaternionic geometry,
  math.DG/0604219, to appear in Tokyo J. Math.

\bibitem{gukovsparks} S. Gukov, J. Sparks: M-Theory on $Spin(7)$
  manifolds, Nucl. Phys. B \textbf{625} (2002), 3--69.

\bibitem{gukospartong} S. Gukov, J. Sparks, D. Tong: Conifold
  transitions and five-brane condensation in $M$-theory on $
    Spin(7)$ manifolds, Class. Quantum Grav. \textbf{20} (2003),
  665--705.

\bibitem{harveylawson} R. Harvey, H.B. Lawson: Calibrated Geometries,
  Acta Math. \textbf{148} (1982), 47--157 

\bibitem{hptt1} E. Heintze, R. Palais, C.-L. Terng, G. Thorbergsson:
  Hyperpolar actions on symmetric spaces, Geometry, topology and
  physics for Raoul Bott, (S.-T. Yau, ed.), International Press,
  Cambridge, (1995)

\bibitem{hitc1} N. Hitchin: Stable forms and special metrics,
  Global differential geometry: the mathematical legacy of Alfred Gray
  (Bilbao, 2000), Contemp. Math. \textbf{288}, 70--89.

\bibitem{hornjohnson} R.A. Horn, C.R. Johnson: \emph{Matrix Algebra},
  Cambridge Univ. Press, 1985.

\bibitem{joyc1} D. Joyce: Compact hypercomplex and quaternionic
  manifolds.  J. Diff. Geom. \textbf{35} (1992), 743--761.

\bibitem{kobswn1} P.Z. Kobak, A.F. Swann: Quaternionic geometry of a
  nilpotent variety, Math. Ann. \textbf{297} (1993), 747--764. 

\bibitem{kollross} A.A. Kollross: A classification of hyperpolar and
  cohomogeneity one actions, Trans. Am. Math. Soc. \textbf{354} (2001)
  571--612.

\bibitem{miyaoka} R. Miyaoka: Bryant-Salamon's $G_{2}$ manifolds and
  the hypersurface geometry, math-ph/0605074

\bibitem{mostert} P.S. Mostert: On a compact Lie group action on
  manifolds, Ann. Math. \textbf{65} (1957), 447--455.

\bibitem{podestaverdiani} F. Podest\`a, L. Verdiani: Positively curved
  $7$-dimensional manifolds, Quart. J. Math. Oxford \textbf{50} (1999), 497--504 

\bibitem{podestaverdiani2} F. Podest\`a, L. Verdiani: Totally geodesic
  orbits of isometries, Ann. Global Anal. Geom., \textbf{16} (1998),
  399--412.

\bibitem{poonsal1} Y.S. Poon, S.M. Salamon: Eight-dimensional
  quaternionic-K\"{a}hler manifolds with positive scalar curvature,
  J. Diff. Geom. \textbf{33} (1991), 363--378.

\bibitem{sal85} S.M. Salamon: Minimal surfaces and symmetric spaces.
  Differential geometry (Santiago de Compostela, 1984), 103--114, Res.
  Notes in Math. 131, Pitman, Boston, MA, 1985.

\bibitem{swann98} A.F. Swann: Homogeneous twistor spaces and nilpotent
  orbits, Math. Ann. \textbf{313} (1999), 161--188. 

\bibitem{uchida} F. Uchida: Classification of compact tranformation
  groups on cohomology complex projective spaces with codimension one
  orbits, Japan J. Math., Vol. 
\textbf{3} (1977), 141--189.

\bibitem{wolf65} J.A. Wolf: Complex Homogeneous contact structures and
  quaternionic symmetric spaces, J. Math. Mech. \textbf{14} (1965),
  1033--1047.

\end{thebibliography}
\end{document}